\documentclass[11pt]{article}
\usepackage{amsmath, amssymb, amsthm, latexsym}
\usepackage{graphicx}
\def \R{{\mathbb R}}

\DeclareMathOperator{\argmin}{argmin}
\usepackage[usenames]{color}
\def \argmin{{\arg\min}}

\usepackage[ruled,vlined]{algorithm2e}

\usepackage{float}
\usepackage{url}

\usepackage{hyperref}

\usepackage{subfloat}
\usepackage{subfig}
\usepackage{overpic}

\usepackage[numbers]{natbib}
\usepackage{authblk}

\begin{document}
\date{}
\title{Extrinsic Bayesian Optimizations  on Manifolds}

\author[1]{Yihao Fang }
\author[2]{Mu Niu }
\author[3]{Pokman Cheung}
\author[1]{Lizhen Lin}

\affil[1]{Department of Applied and Computational Mathematics and Statistics, University of Notre Dame}
\affil[2]{School of Mathematics and Statistics, University of Glasgow}
\affil[3]{London, UK}


\maketitle
\begin{abstract}
    We propose an extrinsic Bayesian optimization (eBO) framework for general optimization problems on manifolds. Bayesian optimization algorithms build a surrogate of the objective function by employing Gaussian processes and  quantify the uncertainty in that surrogate by deriving an acquisition function. This acquisition function represents the probability of improvement based on the kernel of the Gaussian process, which guides the search in the optimization process. The critical challenge for designing Bayesian optimization algorithms on manifolds lies in the difficulty of constructing valid covariance kernels for Gaussian processes on general manifolds. Our approach is to employ extrinsic Gaussian processes by first embedding the manifold onto some higher dimensional Euclidean space via equivariant embeddings and then constructing a valid covariance kernel on the image manifold after the embedding.  This leads to efficient and scalable algorithms for optimization over complex manifolds. Simulation study and real data analysis are carried out to demonstrate the utilities of our eBO framework by applying the eBO to various optimization problems over  manifolds such as the sphere, the Grassmannian, and the manifold of positive definite matrices.
\end{abstract}

\section{Introduction}
   
 Optimization concerns about best decision making which is present in almost aspects of society.  Formally speaking, it aims to optimize some criterion, called the {\it objective function}, over some parameters of variables of interest. In many cases, the variable to optimize possesses certain constraints which should be incorporated  or respected in the optimization process. There is a literature on constrained optimization incorporating linear and nonlinear constraints, including Lagrange-based algorithms and interior points methods.  Our work focuses on an important class of optimization problems with {\it geometric constraints}  in which the parameters or variable to be optimized are assumed to lie on some manifolds, a well-characterized object in differential geometry. In other words, we deal with {\it optimization problems on manifolds}. Optimization on manifolds has abundant applications in modern data science. This is motivated by the systematic collection of modern complex data that take the manifold form. For example,  one may encounter data in the forms of  \emph{positive definite matrices} \citep{dti-ref}, shape objects \citep{kendall77}, subspaces \citep{brian2014,Nishimori2008}, networks \citep{2017eric} and  orthonormal frames \citep{vecdata}. Statistical inference and learning of such data sets often involve optimization problems over manifolds. One of the notable examples is the estimation of the Fr\'echet mean for statistical inference,  which can be cast as an optimization of the Fr\'echet function over manifolds \citep{Frechet1948, bhattacharya2017omnibus, rabimono}. In addition to various examples with data and parameters lying on manifolds, many learning problems in big data analysis with the primary goal of  extracting some lower-dimensions structure, and this lower-dimensional structure is often assumed to  be a manifold. Learning this lower-dimensional structure often requires solving optimization problems over manifolds such as the Grassmannian.

The critical challenge for solving \emph{optimization problems over manifolds} lies in how to appropriately incorporate the underlying geometry of manifolds  for optimization.  Although there has been a fast development in optimization (over Euclidean spaces in general), and it is an extremely active ongoing research area,  there is a tremendous challenge for extending theories and algorithms developed in optimization over Euclidean spaces to  manifolds. The optimization approach on manifolds is superior to performing free Euclidean optimization and projecting the parameters back onto the search space after each iteration, such as in the projected gradient descent method. It has been shown to outperform standard algorithms for many problems. Some of those algorithms, such as Newton method \citep{2014arXiv1407.5965S, ring12},  conjugate gradient descent algorithm \citep{ Nishimori2008}, steepest descent \citep{AbsMahSep2010} and trusted region method, \citep{Absil2007}  have recently been extended to the Riemannian manifold from the Euclidean space, and most of the methods require the knowledge of the gradients. Further successful applications of optimization on manifolds include efficient parallel algorithm \citep{saparbayeva2018communication}, accelerated algorithm \citep{lin2020accelerated}, and natural gradient algorithm \citep{izadi2020optimization}. However, in many cases, analytical or simple forms of gradient information are unavailable. In other cases, the evaluations and calculations of the gradient or the Hessian (second form for the case of the manifold) can be expensive. There is a lack of gradient-free methods for optimization problems when the  gradient information is not available or expensive to obtain, or the objective function is expensive to evaluate. In these scenarios, gradient-free methods can be appealing alternatives.  Many gradient-free methods are proposed for optimization problems in the Euclidean space, especially the state-of-the-art Bayesian optimization method, which has emerged as a very powerful tool in machine learning for tuning both learning parameters and hyperparameters \citep{NIPS2012_4522}. Such algorithms outperform many other global optimization algorithms \citep{Jones:2001}.  Bayesian optimization originated with the work of Kushner and Mockus \cite{kushner1964new,movckus1975bayesian}. It received considerably more attention after Jones' work on the Efficient Global Optimization (EGO) algorithm \cite{jones1998efficient}. Following that work, innovations developed in the same literature include multi-fidelity optimization \cite{huang2006sequential}, multi-objective optimization \cite{keane2006statistical}, and a study of convergence rates \cite{calvin1997average}. The observation by Snoek \cite{snoek2012practical} that Bayesian optimization is useful for training deep neural networks sparked significant interest in machine learning. Within this trend, Bayesian optimization has also been used to choose laboratory experiments in materials and drug design \cite{packwood2017bayesian}, in the calibration of environmental models \cite{shoemaker2007watershed}, and in reinforcement learning \cite{brochu2010tutorial}.

However, the overwhelming literature mainly focuses on the input domain in a Euclidean space. In recent years, with the surging collection of complex data, it's common for the input domain doesn't have such a simple form. For instance, the inputs might be restricted to a non-Euclidean manifold. \cite{lin2019extrinsic} proposed a general extrinsic framework for GP modeling on manifolds, which depends on the embedding of the manifold in the Euclidean space and the construction of extrinsic kernels for GPs on their images. \cite{jaquier2020high} learned a nested-manifolds embedding and a representation of the objective function in the latent space using Bayesian optimization in low-dimensional latent spaces. \cite{borovitskiy2020matern} exploited the geometry of non-Euclidean parameter spaces arising in robotics by using Bayesian optimization to properly measure similarities in the parameter space through the Gaussian process.
Recently, \cite{jaquier2022geometry} implemented Riemannian Matern kernels on manifolds to place the Gaussian process prior and applied these kernels on the geometry-aware Bayesian optimization for a variety of robotic applications, including orientation control, manipulability optimization, and motion planning.

Motivated by the success of Bayesian optimization algorithms, in this work, we propose to develop the eBO framework on manifolds. In particular, we employ an extrinsic class of Gaussian processes on manifolds as the surrogate model for the objective function and quantify the uncertainty in that surrogate. An acquisition function can also be defined from this surrogate to decide where to sample \cite{frazier2018tutorial}. The algorithms are demonstrated and applied to both simulated and real data sets to various optimization problems over different manifolds, including spheres, positive definite matrices, and the Grassmannian.

We organize our work as follows. In Section 2, we provide a streamlined introduction to Bayesian optimization on manifolds employing the extrinsic Gaussian process (eGP) on manifolds. In Section 3, we present a concrete illustration of our eBO algorithm in the extrinsic setting in terms of optimizing the acquisition function from eGP. In Sec. 4, we demonstrate numerical experiments on optimization problems over different manifolds, including  spheres, positive definite matrices, and the Grassmann manifold, and showcase the performance of our approach.

\section{ Bayesian optimization (BO) on manifolds}

Let $f(x)$ be an objective function defined over some manifold $M$.  We are interested in solving
\begin{align}
\mu= \argmin_{x\in M}f(x).
\end{align}
Typically, the function $f$ lacks a known special structure such as concavity or linearity that would make it easy to optimize using techniques that leverage such structure to improve efficiency. Furthermore, we assume one only needs to be able to evaluate $f(x)$ without having to know first or second order information such as gradient  when evaluating $f$. We refer to problems having this property as ``derivative-free". Lastly, we seek the global optimizer if it exists. The principal ideas behind Bayesian optimization are to  build a probabilistic model for the objective function by imposing a \emph{Gaussian process prior}, and this probabilistic model will be used to guide to where in $M$ the function is next evaluated. The  posterior predictive distribution will then be computed. Instead of optimizing the usually expensive objective function, a cheap proxy function is often optimized, which will determine the next point to evaluate. One of the inherent difficulties lies in constructing \emph{valid Gaussian processes on manifolds} that will be utilized for Bayesian optimization on manifolds.

To be more specific, as above, let $f(x)$ be the objective function  on $M$ for which we omit the  potential dependence on some data, for now, the goal is to find 
the minimum:  $\mu = \arg \min_{x\in M} f (x)$. 
Let $f(x)\sim GP(\nu(\cdot), R(\cdot, \cdot))$, $x\in M$,  be a Gaussian process (GP) on the manifold $M$ with mean function $\nu(x)$ and covariance kernel $R(\cdot,\cdot)$. Then we evaluate $f(x)$  at a finite number of points on the manifold following a multivariate Gaussian distribution, that is, 
\begin{align*}
(f(x_1),\ldots, f(x_n))\;\sim\; N\left((\nu(x_1),\ldots, \nu(x_n)), \Sigma\right), \\
\;\Sigma_{ij}=cov\left(f(x_i), f(x_j)\right)=R(x_i,x_j).
 \end{align*}
Here  $R(\cdot,\cdot):M\times M\rightarrow \R$  is  a  covariance kernel defined on the manifold, which is a \emph{positive semi-definite kernel} on  $M$. It states that, for any sequence $(a_1,\ldots, a_n)\in \R^n$, $\sum_{i=1}^n\sum_{j=1}^n a_ia_jR(x_i,x_j)\geq 0.$

After some evaluations, the GP gives us closed-form marginal means and variances. We denote the predictive means and variance as $\nu(x, \mathcal D)$ and $\sigma^2(x, \mathcal D)$, where $\mathcal D=\{x_1,\ldots, x_n\}$.  The acquisition function, which we denote by $a : M\rightarrow \R^+$, determines which point in $M$ should be evaluated next via a proxy optimization 
\begin{equation}
x_{next} = \arg\max_{x\in M} a(x).
\end{equation}There are several popular choices of acquisition functions that are proposed for Bayesian optimization in the Euclidean space \cite{frazier2018tutorial}. The analogous version can be generalized to manifolds.  Under the Gaussian process prior, these functions depend on the model solely through its predictive mean function  $\nu(x, \mathcal D)$ and  variance $\sigma^2(x, \mathcal D)$.  In the preceding, we denote the best current value as  
\begin{equation}x_{best} = \argmin_{x_n} f(x_n). 
\end{equation} We denote the cumulative distribution function of the standard normal as $\Phi(\cdot)$.  We adopt one of the acquisition functions as
\begin{align}
\label{acquisition}
a(x)=\Phi(r(x)), \;\text{where}\; r(x)=\dfrac{f(x_{best})-\nu(x, \mathcal D)}{\sigma(x, \mathcal D)},
\end{align}
which represents \emph{probability of improvement}.
Other acquisition functions, such as knowledge gradient and entropy search on manifolds, could also be developed, but here we focus on probability improvement. As seen above, a key component of Bayesian optimization methods on manifolds is utilizing Gaussian processes on manifold $M$, which are non-trivial to construct. The critical challenge for constructing the Gaussian process is constructing \emph{valid covariance kernel} on manifolds. In the following sections, we will utilize extrinsic Gaussian processes on manifolds by employing extrinsic covariance kernels via embeddings \citep{lin2019extrinsic} on the specific manifold.

\section{Extrinsic Bayesian Optimization on manifolds}
In this section, we will present the extrinsic Bayesian optimization algorithm on manifolds with the help of the illustrated eGP on manifolds. The key component of the extrinsic framework is the embedding map. Let $J: M\rightarrow\R^D$ be an embedding of $M$ into some higher dimensional Euclidean space $\mathbb{R}^D$ ($D\geq d$) and denote the image of the embedding as $\tilde{M} = J(M)$. By definition of an embedding, $J$ is a smooth map such that it's differential  $dJ: T_xM\rightarrow T_{J(x)}\R^D$ at each point $x\in M$ is an injective map from its tangent space $T_xM$  to  $T_{J(x)}\R^D$, and $J$ is a homeomorphism between $M$ and its image $\tilde{M}$. As mentioned above, we use eGP as the prior distribution of the objective function, and the main algorithm is described in algorithm\ref{ebo}. Let $f(x)\sim eGP(m(\cdot), R(\cdot, \cdot))$ with prior mean $m(x)$.  $R(\cdot, \cdot)$ is an covariance kernel on $M$ induced by the embedding $J$ as:
\begin{align}
R(x,z)=\tilde R(J(x), J(z)),
\end{align}
where $\tilde R(\cdot, \cdot)$ is a valid covariance kernel on $\mathbb R^D$.

As above, let the predictive means and variance be $\nu(x, \mathcal D)$ and $\sigma^2(x, \mathcal D)$, where $\mathcal D$ denotes the data. The acquisition function defined by equation \eqref{acquisition} represents the probability of improvement. To find out the next point to evaluate, we need to maximize the acquisition question. Since we use the zero mean eGP where $m(\cdot)$ equals 0,  the expression of  $\nu(x, \mathcal D)$ and $\sigma(x, \mathcal D)$ are given as follows:
\begin{align}
\nu(x, \mathcal D)&=  K(x_D,x)^T (K(x_D,x_D)+ \sigma_n^2 I )^{-1}y ,\\
\sigma(x, \mathcal D)&= K(x,x) - k(x_D,x)^T ( K(x_D,x_D) + \sigma_n^2 I )^{-1} k(x_D,x).
\end{align}
Noting that $a(x)$ depends on $x$ through $J(x)$, let $a(x)=\tilde a(\tilde x)$ where $\tilde x=J(x)\in \tilde M$. Therefore, the optimization of $a(x)$ over $x$ is \emph{equivalent to optimization of $\tilde a (\tilde x)$ over $\tilde M$.} Let $\tilde x^*$ be the optimizer of $\tilde a(x)$ over $\tilde M$, i.e.,
\begin{align}
\label{eq-exoptimizer}
\tilde x^*=\argmin_{\tilde x\in \tilde M}\tilde a(\tilde x).
\end{align}
Then one has
\begin{align}
x^*=J^{-1}(\tilde x^*),
\end{align}
where 
\begin{align}
x^*=\argmin_{x\in M}a(x).
\end{align}

The key is to solve \eqref{eq-exoptimizer}.  We consider the gradient descent or Newton's method \emph{over the submanifolds $\tilde M$} since the gradient is easy to obtain on $\tilde M$. Let $grad\;\tilde a (\tilde x)$ in the Euclidean space, then the gradient of $\tilde a(\tilde x)$ over the submanifold $\tilde M$ is given by
\begin{align}
\mathcal Pgrad\;\tilde a (\tilde x),
\end{align}
where $\mathcal P$ is the projection map $P$ from $T_{\tilde x}\R^D$ onto the tangent space $T_{\tilde x}\tilde M$. We propose the following gradient descent algorithm \ref{gd} for finding $\tilde x^*$ of $\tilde a(\tilde x)$ over $\tilde M$ by following steps

\begin{itemize}
\item Let $\tilde x_0\in \tilde M$ be an initial point. 
\item  $\tilde x_{t+1}=\exp_{\tilde x_t}\lambda \mathcal Pgrad\;\tilde a (\tilde x)\mid_{\tilde x_{t}}$
where $\exp$ is the exponential map on the submanifold. This exponential map is with respect to the metric by restricting the Euclidean metric onto the image manifold or submanifold. 
\end{itemize}

\begin{algorithm}
	\caption{Gradient algorithms for optimization $\tilde a (\tilde x)$ along the submanifold $\tilde M$}
\label{gd}
	\For{$t=0, 1, \dots, T$}{
Let $\tilde x_0\in \tilde M$ be an initial point. \\
$\tilde x_{t+1}=\exp_{\tilde x_t}\lambda \mathcal Pgrad\;\tilde a (\tilde x)\mid_{\tilde x_{t}}$,\\
}
where $\exp$ is the exponential map on the submanifold. 
\end{algorithm}

\begin{algorithm}
		\caption{Extrinsic Bayesian optimization on manifolds}
\label{ebo}
	Initialize $x_1,\ldots, x_k\in \tilde M$ ($k\geq 1$);\\
	Let $x_{best}^0=\argmin\{f(x_1),\ldots, f(x_k)\}$ and $\mathcal D=\{x_1,\ldots, x_k\}$;\\
	\For{$s=0, 1, \dots, T-1$}{
$x_{next}=J^{-1}(\arg\max	\tilde a (\tilde x))$;\\
$x_{best}^{s+1}=\argmin\{(f(x_1),\ldots, f(x_k), f(x_{next})\}$;\\

Update $\mathcal D=\{\mathcal D, x_{next}\}$
		}
	Return $x_{best}^T$
\end{algorithm}

\section{Examples and Applications}

To illustrate the broad applications of our eBO framework, we utilize a large class of examples with data domains on different manifolds, including spheres, positive definite matrices, and the Grassmann manifold. We construct the extrinsic kernels for eGPs based on the corresponding embedding map on the specific manifold from \cite{lin2019extrinsic}. In Section 4.1, the Fr\'echet mean estimation is carried out with data distributed on the sphere. In Section 4.2, the eBO method is applied to the matrix approximation problem on Grassmannians. Lastly, Section 4.3 considers a  regression problem on the manifold of positive definite matrices, which has essential applications in neuroimaging. We use eBO method to show the difference between healthy samples and HIV+ samples.

\subsection{Estimation of Fr\'echet means}

We will first apply the eBO method to the estimation of sample Fr\'echet means on the sphere. Let $x_1,\ldots, x_n$ be $n$ points on the manifold, such as the sphere.  The sample Fr\'echet function is defined as
 \begin{align}
\label{eq-empiricalfrechet}
f_n(x)=\frac{1}{n}\sum_{i=1}^n\rho^2(x,x_i),
\end{align}
and the minimization of $f_n(x)$ leads to the estimation of sample Fr\'echet mean $\mu_n$, where
\begin{align}
\label{frechet-sphere}
\mu_n=\arg\min_{x\in M} \frac{1}{n}\sum_{i=1}^n\rho^2(x,x_i).
\end{align}

Evaluation of  $f_n(x)$ at the data point $x_1,\ldots, x_n$ will be our data points for guiding the optimization to find the Fr\'echet mean. The choice of distance $\rho$ leads to different notions of means on manifolds. In the extrinsic setting, we can define an extrinsic mean via some embedding of the manifold onto some Euclidean space.  Specifically, let $J: M \rightarrow \R^D$ be some embedding of the manifold onto some higher-dimensional Euclidean space $\R^D$. Then one can define the extrinsic distance as 
\begin{align}
\rho(x,z)=\|J(x)-J(z)\|,
\end{align}
where $\|\cdot\|$ is the Euclidean distance on $\R^D$. This leads to the  extrinsic Fr\'echet mean. Fortunately, the extrinsic mean has a close form expression, namely $\mu_n=\mathcal P_{\tilde M}\left(\frac{\sum_{i=1}^n J(x_i)}{n}\right)$, where $\mathcal P$ stands for the projection map onto the image $\tilde M=J(M)$.

\begin{figure}[htbp]
\centering
\includegraphics[width=10cm, height = 8cm]{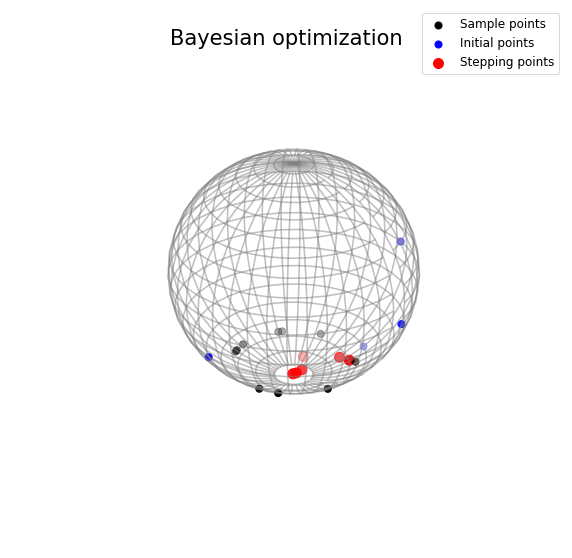}
\caption{The iteration steps of the eBO method on the sphere $S^2$. The data points $x_1,...,x_n$ are plotted as those black points on the same latitude of the sphere. The true extrinsic mean is the south pole $(0,0,-1)$. We start with those random blue points on the sphere. The outputs of iterations in our algorithm \ref{ebo} are marked as red points on the sphere, converging to the ground truth. \label{fig:sphere1}}
\end{figure}

\begin{figure}[htbp]
\centering
\includegraphics[width=12cm, height = 6cm]{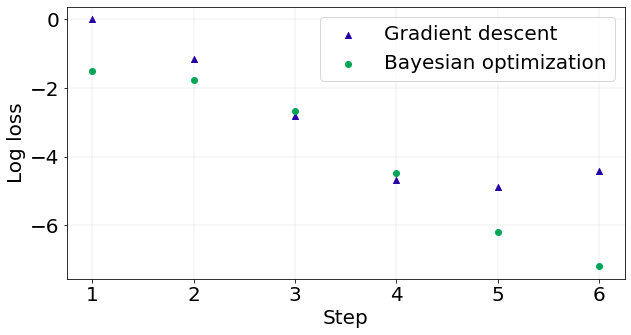}
\caption{We compare the eBO method with the gradient descent (GD) method on the sphere by calculating the log L2 error from the true extrinsic mean. Although the eBO method converges slower than the GD method in the first four steps, by adding those stepping points into the eGP, it achieves better numerical precision with more steps.    \label{fig:sphere2}}
\end{figure}

In our simulation study, we estimate the extrinsic mean on the sphere $S^2$ and evaluate the performance of our eBO method by comparing it to the gradient descent (GD) algorithm in terms of convergence to the ground truth (the true minimizer given above). The embedding map $J$ is the identity map $I: \mathbb S^2 \rightarrow \mathbb R^3$. As shown in Fig \ref{fig:sphere1}, we simulate some data on the sphere, and the goal is to find the extrinsic mean by solving equation \eqref{frechet-sphere}. Specifically, those $x_i$ $(i=1, \dots, n)$ points in black are sampled on the circle of latitude, which leads to the south pole $(0,0,-1)$ as the true extrinsic Fr\'echet mean. In terms of our BO algorithm \ref{ebo}, we first sample those blue points randomly on the sphere as initial points and initialize the covariance matrix by evaluating the covariance kernel at these points. Then, we randomly select the direction on the sphere to minimize the acquisition function based on the eGP.  In each iteration, we mark the minimizer as the stepping point in  Figure \ref{fig:sphere1} and add it to the data for the next iteration. Not surprisingly, those stepping points in red converge to the ground truth (south pole) after a few steps. We also compare our eBO method with the gradient descent (GD) method on the sphere under the same initialization. Illustrated in Fig \ref{fig:sphere2}, although our BO method converges slightly slower than GD at the first four steps, the BO method achieves better accuracy with a few more steps. It confirms the quick convergence and high accuracy as the advantages of the eBO method.

\subsection{Invariant subspace approximation on Grassmann manifolds}

We investigate two important manifolds, the Stiefel manifolds and the Grassmann manifolds (Grassmannians). The Stiefel manifold is the collection of $p$ orthonormal frames in $\R^n$, that is, $Stiefel(p,n)= \{X\in \mathbb R^{n*p}: X^TX = I_p \}$. Moreover, the Grassmann manifold is the space of all the subspaces of a fixed dimension $p$ whose basis elements are $p$ orthonormal unit vectors in $\mathbb R^n$, which is $Grass(p,n) = \{span(X): X\in \mathbb R^{n*p}, X^TX = I_p \}$. Those two manifolds are closely related. The key difference between a point on the Grassmannian and a point on the Stiefel manifold is that the ordering of the $p$ orthonormal vectors in $\mathbb R^n$ does not matter for the former. In other words, the Grassmannian could be viewed as the quotient space of the Stiefel manifold modulo $O(p)$, the $p*p$ orthogonal group. That is,  $Grass(p,n) = Stiefel(p,n)/O(p)$.

We consider the matrix approximation problem on the Grassmannian subspace manifold and apply the eBO method to solve it. Given a full rank matrix $F \in \mathbb R^{n*m}$, without loss of generality, we assume $n<=m$ and $rank(F)=n$, the goal is to approximate this matrix $F$ in the subspace $Grass(p,n)$ with $p<m$. From any matrix $X\in Grass(p,n)$, we approximate the original matrix $F$ from the algorithm \ref{grass}.

\begin{algorithm}
		\caption{Matrix approximation on the Grassmann manifold}
\label{grass}
	Initialize the orthogonal matrix $X\in Grass(p,n)$;\\
	Initialize the matrix $W\in \mathbb R^{p*m}$ and calculate the column $W_j \in \mathbb R^{r*1}$ below;\\
	\For{$j=1, \dots, m$}{
 $W_{j}= (X^TX)^{-1} X^TF_j$, where $F_j$ is the $j$th column of $F$;\\
		}
	Return $XW$ as the approximation.
\end{algorithm}

Here we consider the approximation error in Frobenius norm where $W$ depends on $X$:
\begin{align}
\label{Grass}
X = \arg\min_{X\in Grass(p,n) } \|XW-F \|_F.
\end{align}

\begin{figure}[htbp]
\centering
\includegraphics[width=12cm, height = 6cm]{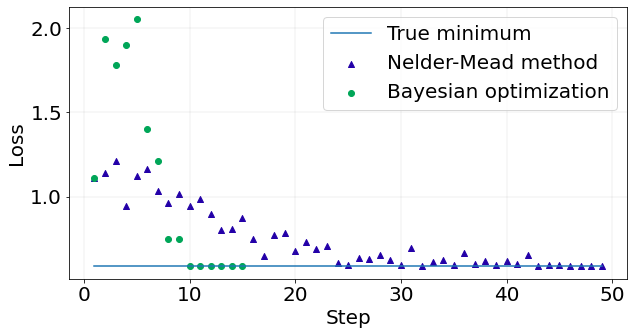}
\caption{We compare the eBO method with the Nelder-Mead method on the matrix approximation problem on the Grassmannian. The minimum of equation\ref{Grass} is around $0.5578$ and plotted as the blue line in the figure. It can't achieve zero loss due to the low dimension constraint. The Nelder-Mead method achieves the minimal subspace around 25 steps and becomes stable after 40 steps. On the other hand, the eBO method converges to the minimum after 10 steps, much faster than the Nelder-Mead method, showing the quick convergence to the optimal solution.   \label{fig:Grassmannian}}
\end{figure}

 In our simulation, we consider the matrix $F \in \mathbb R^{3*6} $ with $rank(F)=3$ and $p=2$. As we know, this matrix approximation problem achieves its minimal when $X \in span(U_1,U_2)$ where $U$ is the left matrix in the Singular Value Decomposition (SVD) and $U_j$ denotes its $j$th column. Let $\hat{U} = [U_1,U_2]$ contain the first two columns of the SVD part $U$. To apply the eBO method, we first sample some matrix following $X_i = \hat{U} + \frac{i*(-1)^{i}}{2} $, $i=1,...,6$ as initial points. Similar to the sphere case, we have to map the minimizer of the acquisition function back to the  manifold via the inverse embedding $J^{-1}$. In $Grass(2,3)$, for a given matrix $X \in \mathbb R^{3*2}$, $J^{-1}$ indicates us to apply the SVD decomposition to $X^TX$ then keep the first two columns. Since the closed form of the loss gradient is untouchable, we compare BO method with the Nelder-Mead method. As shown in Figure \ref{fig:Grassmannian}, our eBO method converges to the ground truth in a few steps, faster than the Nelder-Mead method.

\subsection{Positive definite matrices}

Lastly, we apply the eBO method to a regression problem with the response on positive definite matrices. Let $(z_1, y_1),\dots, (z_n, y_n)$ be $n$ points from a regression model in which $z_i\in \R^d$ and $y\in M=SPD(p)$, where $SPD(p)$ stands for $p$ by $p$ positive semi-definite matrices. We are interested in modeling the regression map between $z$ and $y$. 

Let $K(\cdot,\cdot): \R^d\times \R^d \rightarrow \R$ be a kernel function defined on  covariate space $\R^d$. For example, one can take the standard Gaussian kernel. Let $g(z)\in M$ be the regression map evaluated at a covariate level $z$. We propose to estimate $g$ as 
\begin{align}
\label{kernel-reg}
g(z)=\arg\min_{y\in M}\sum_{i=1}^n\frac{1}{h}K_h(z, z_i) \rho^2(y, y_i),
\end{align}
where $h$ is the bandwidth parameter of the kernel function $K$.  We take $\rho$ to be some intrinsic distance on the SPD. For the given input or covariate $z$, we denote $f_n(y)=\sum_{i=1}^n\frac{1}{h}K_h(z, z_i) \rho^2(y, y_i))$ and apply our Bayesian optimization technique in minimizing  the objective function  $f_n(y)$  whose minimizer turns out to be the estimate of our regression function for $z$. 

In our analysis, we focus on the case of $p=3$, which has important applications in diffusion tensor imaging (DTI), which is designed to measure the diffusion of water molecules in the brain. In more detail, diffusion represents the directional along the white matter tracks or fibers, corresponding to structural connections between brain regions along where brain activity and communications happen. DTI data are collected routinely in human brain studies. People are interested in using DTI to build predictive models of cognitive traits, and neuropsychiatric disorders \cite{lin2019extrinsic}. The diffusions are characterized in terms of diffusion matrices, represented by 3 by 3 positive definite matrices. The data set consists of $46$ subjects with $28$ $HIV+$ subjects and $18$ healthy controls. Those 3 by 3 diffusion matrices were extracted along one atlas fiber tract of the splenium of the corpus callosum. All the DTI data were registered in the same atlas space based on arc lengths, with $75$ tensors obtained along the fiber tract of each subject, which has been studied in \cite{Yuan2012,lin2019extrinsic, ereg} in a regression and kernel regression setting. Here we consider the arc length tensor as the covariant $z$ and the diffusion matrix as $y$ in our framework. We aim to show the difference between the positive (HIV+) samples and control samples based on the weighted Fr\'echet mean at each arc length. 

\begin{figure}[htbp]
\centering
\includegraphics[width=12cm, height = 6cm]{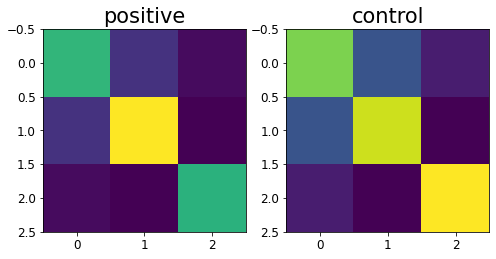}
\caption{Estimated diffusion tensors (DTI) at first arc length are shown as $3*3$ matrix from the healthy group and the HIV+ group. The horizontal and vertical axis denotes the rows and columns of the matrix. Entries inside the matrix are represented in different colors. Based on the color bar, differences between healthy and HIV+ groups could be observed, especially on the diagonal elements.    \label{fig:psd}}
\end{figure}

We apply our eBO methods to build the kernel regression model with the positive group and the control group separately. That is, for each $z$ in the arc lengths (locations), from \ref{kernel-reg}, we obtain estimated diffusion tensor $g_{1}(z)$ as the weighted Fr\'echet mean from the positive samples, and $g_{2}(z)$ as the weighted Fr\'echet mean from the control samples as well. Furthermore, the DTIs at the first arc length $g_{1}(z_1)$ and $g_{2}(z_1)$ from different groups are shown in figure \ref{fig:psd}, where we observe the value difference between the positive sample and the control sample. Consequently, we carry out the two sample test \cite{bhattacharya2017omnibus} in terms of the estimated diffusion tensors from all $75$ arc lengths and yield an extremely small $p-value$ equal to $O(10^{-5})$. It indicates the difference between those two groups based on our BO estimations.

\section{Discussion and conclusion}

We propose a general extrinsic framework for Bayesian optimizations on manifolds. These algorithms are based on the Gaussian process's acquisition function on manifolds. Applications are presented by applying the eBO method to various optimization and regression problems with data on different manifolds, including spheres, Stiefel/Grassmann manifolds, and the spaces of positive definite matrices. As a gradient-free approach, the eBO method shows advantages compared to the gradient descend method and Nelder-Mead method in our simulation study. In future work, we will investigate the intrinsic Bayesian optimizations on manifolds based on the intrinsic Gaussian processes such as the ones based on heat kernels \cite{Niu2018IntrinsicGP}.

\section*{Acknowledgement}
We acknowledge the generous support of NSF grants DMS CAREER 1654579 and DMS 2113642.

\bibliographystyle{abbrv}
\bibliography{reference}

\begin{thebibliography}{10}

\bibitem{Absil2007}
P.-A. Absil, C.~Baker, and K.~Gallivan.
\newblock Trust-region methods on riemannian manifolds.
\newblock {\em Foundations of Computational Mathematics}, 7(3):303--330, Jul
  2007.

\bibitem{AbsMahSep2010}
P.-A. Absil, R.~Mahony, and R.~Sepulchre.
\newblock Optimization on manifolds: methods and applications.
\newblock In M.~Diehl, F.~Glineur, and W.~Michiels, editors, {\em Recent Trends
  in Optimization and its Applications in Engineering}. Springer-Verlag, 2010.

\bibitem{dti-ref}
A.~Alexander, J.~Lee, M.~Lazar, and A.~Field.
\newblock Diffusion tensor imaging of the brain.
\newblock {\em Neurotherapeutics}, 4(3):316--329, 2007.

\bibitem{rabimono}
A.~Bhattacharya and R.~Bhattacharya.
\newblock {\em Nonparametric Inference on Manifolds: with Applications to Shape
  Spaces}.
\newblock Cambridge University Press, 2012.
\newblock IMS monographs \#2.

\bibitem{bhattacharya2017omnibus}
R.~Bhattacharya and L.~Lin.
\newblock Omnibus clts for fr{\'e}chet means and nonparametric inference on
  non-euclidean spaces.
\newblock {\em Proceedings of the American Mathematical Society},
  145(1):413--428, 2017.

\bibitem{borovitskiy2020matern}
V.~Borovitskiy, A.~Terenin, P.~Mostowsky, et~al.
\newblock Mat{\'e}rn gaussian processes on riemannian manifolds.
\newblock {\em Advances in Neural Information Processing Systems},
  33:12426--12437, 2020.

\bibitem{brochu2010tutorial}
E.~Brochu, V.~M. Cora, and N.~De~Freitas.
\newblock A tutorial on bayesian optimization of expensive cost functions, with
  application to active user modeling and hierarchical reinforcement learning.
\newblock {\em arXiv preprint arXiv:1012.2599}, 2010.

\bibitem{calvin1997average}
J.~M. Calvin.
\newblock Average performance of a class of adaptive algorithms for global
  optimization.
\newblock {\em The Annals of Applied Probability}, pages 711--730, 1997.

\bibitem{vecdata}
T.~D. Downs, J.~Liebman, and W.~Mackay.
\newblock Statistical methods for vectorcardiogram orientations.
\newblock In R.~H. I.~Hoffman and E.~E.~Glassman, editors, {\em
  Vectorcardiography 2: Proc. XIth Intn. Symp. Vectorcardiography}, pages
  216--222. North-Holland, Amsterdam., 1971.

\bibitem{frazier2018tutorial}
P.~I. Frazier.
\newblock A tutorial on bayesian optimization.
\newblock {\em arXiv preprint arXiv:1807.02811}, 2018.

\bibitem{Frechet1948}
M.~Fr\'{e}chet.
\newblock Les \'{e}l\'{e}ments al\'{e}atoires de nature quelconque dans un
  espace distanci\'{e}.
\newblock {\em Annales de L'{I}nstitut {H}enri {P}oincar\'{e}}, 10(4):215--310,
  1948.

\bibitem{huang2006sequential}
D.~Huang, T.~T. Allen, W.~I. Notz, and R.~A. Miller.
\newblock Sequential kriging optimization using multiple-fidelity evaluations.
\newblock {\em Structural and Multidisciplinary Optimization}, 32(5):369--382,
  2006.

\bibitem{izadi2020optimization}
M.~R. Izadi, Y.~Fang, R.~Stevenson, and L.~Lin.
\newblock Optimization of graph neural networks with natural gradient descent.
\newblock In {\em 2020 IEEE international conference on big data (big data)},
  pages 171--179. IEEE, 2020.

\bibitem{jaquier2022geometry}
N.~Jaquier, V.~Borovitskiy, A.~Smolensky, A.~Terenin, T.~Asfour, and L.~Rozo.
\newblock Geometry-aware bayesian optimization in robotics using riemannian
  mat{\'e}rn kernels.
\newblock In {\em Conference on Robot Learning}, pages 794--805. PMLR, 2022.

\bibitem{jaquier2020high}
N.~Jaquier and L.~Rozo.
\newblock High-dimensional bayesian optimization via nested riemannian
  manifolds.
\newblock {\em Advances in Neural Information Processing Systems},
  33:20939--20951, 2020.

\bibitem{Jones:2001}
D.~R. Jones.
\newblock A taxonomy of global optimization methods based on response surfaces.
\newblock {\em J. of Global Optimization}, 21(4):345--383, Dec. 2001.

\bibitem{jones1998efficient}
D.~R. Jones, M.~Schonlau, and W.~J. Welch.
\newblock Efficient global optimization of expensive black-box functions.
\newblock {\em Journal of Global optimization}, 13(4):455--492, 1998.

\bibitem{keane2006statistical}
A.~J. Keane.
\newblock Statistical improvement criteria for use in multiobjective design
  optimization.
\newblock {\em AIAA journal}, 44(4):879--891, 2006.

\bibitem{kendall77}
D.~G. Kendall.
\newblock The diffusion of shape.
\newblock {\em Adv. Appl. Probab.}, 9:428--430, 1977.

\bibitem{2017eric}
E.~{Kolaczyk}, L.~{Lin}, S.~{Rosenberg}, and J.~{Walters}.
\newblock {Averages of Unlabeled Networks: Geometric Characterization and
  Asymptotic Behavior}.
\newblock {\em ArXiv e-prints}, Sept. 2017.

\bibitem{kushner1964new}
H.~J. Kushner.
\newblock A new method of locating the maximum point of an arbitrary multipeak
  curve in the presence of noise.
\newblock {\em Journal of Fluids Engineering}, 1964.

\bibitem{lin2019extrinsic}
L.~Lin, N.~Mu, P.~Cheung, and D.~Dunson.
\newblock Extrinsic gaussian processes for regression and classification on
  manifolds.
\newblock {\em Bayesian Analysis}, 14(3):887--906, 2019.

\bibitem{lin2020accelerated}
L.~Lin, B.~Saparbayeva, M.~M. Zhang, and D.~B. Dunson.
\newblock Accelerated algorithms for convex and non-convex optimization on
  manifolds.
\newblock {\em arXiv preprint arXiv:2010.08908}, 2020.

\bibitem{ereg}
L.~Lin, B.~S. Thomas, H.~Zhu, and D.~B. Dunson.
\newblock Extrinsic local regression on manifold-valued data.
\newblock {\em Journal of the American Statistical Association},
  112(519):1261--1273, 2017.

\bibitem{movckus1975bayesian}
J.~Mo{\v{c}}kus.
\newblock On bayesian methods for seeking the extremum.
\newblock In {\em Optimization techniques IFIP technical conference}, pages
  400--404. Springer, 1975.

\bibitem{Nishimori2008}
Y.~Nishimori, S.~Akaho, and M.~D. Plumbley.
\newblock {\em Natural Conjugate Gradient on Complex Flag Manifolds for Complex
  Independent Subspace Analysis}, pages 165--174.
\newblock Springer Berlin Heidelberg, Berlin, Heidelberg, 2008.

\bibitem{Niu2018IntrinsicGP}
M.~Niu, P.~Cheung, L.~Lin, Z.~Dai, N.~D. Lawrence, and D.~B. Dunson.
\newblock Intrinsic gaussian processes on complex constrained domains.
\newblock {\em Journal of the Royal Statistical Society: Series B (Statistical
  Methodology)}, 81:603--627, 2019.

\bibitem{packwood2017bayesian}
D.~Packwood.
\newblock {\em Bayesian optimization for materials science}.
\newblock Springer, 2017.

\bibitem{ring12}
W.~Ring and B.~Wirth.
\newblock Optimization methods on riemannian manifolds and their application to
  shape space.
\newblock {\em SIAM Journal on Optimization}, 22(2):596--627, 2012.

\bibitem{saparbayeva2018communication}
B.~Saparbayeva, M.~Zhang, and L.~Lin.
\newblock Communication efficient parallel algorithms for optimization on
  manifolds.
\newblock {\em Advances in Neural Information Processing Systems}, 31, 2018.

\bibitem{shoemaker2007watershed}
C.~A. Shoemaker, R.~G. Regis, and R.~C. Fleming.
\newblock Watershed calibration using multistart local optimization and
  evolutionary optimization with radial basis function approximation.
\newblock {\em Hydrological sciences journal}, 52(3):450--465, 2007.

\bibitem{2014arXiv1407.5965S}
S.~T. {Smith}.
\newblock {Optimization Techniques on Riemannian Manifolds}.
\newblock {\em ArXiv e-prints}, July 2014.

\bibitem{NIPS2012_4522}
J.~Snoek, H.~Larochelle, and R.~P. Adams.
\newblock Practical bayesian optimization of machine learning algorithms.
\newblock In F.~Pereira, C.~J.~C. Burges, L.~Bottou, and K.~Q. Weinberger,
  editors, {\em Advances in Neural Information Processing Systems 25}, pages
  2951--2959. Curran Associates, Inc., 2012.

\bibitem{snoek2012practical}
J.~Snoek, H.~Larochelle, and R.~P. Adams.
\newblock Practical bayesian optimization of machine learning algorithms.
\newblock {\em Advances in neural information processing systems}, 25, 2012.

\bibitem{brian2014}
B.~{St.~Thomas}, L.~{Lin}, L.-H. {Lim}, and S.~{Mukherjee}.
\newblock {Learning subspaces of different dimension}.
\newblock {\em ArXiv e-prints}, 1404.6841, Apr. 2014.

\bibitem{Yuan2012}
Y.~Yuan, H.~Zhu, W.~Lin, and J.~S. Marron.
\newblock Local polynomial regression for symmetric positive definite matrices.
\newblock {\em Journal of the Royal Statistical Society: Series B (Statistical
  Methodology)}, pages no--no, 2012.

\end{thebibliography}
\end{document}